\newcount\fontset
  \fontset=1
  \def\dualfont#1#2#3{\font#1=\ifnum\fontset=1 #2\else#3\fi}

  \dualfont\bbfive{bbm5}{cmbx5}
  \dualfont\bbseven{bbm7}{cmbx7}
  \dualfont\bbten{bbm10}{cmbx10}
  \font \eightbf = cmbx8
  \font \eighti = cmmi8 \skewchar \eighti = '177
  \font \eightit = cmti8
  \font \eightrm = cmr8
  \font \eightsl = cmsl8
  \font \eightsy = cmsy8 \skewchar \eightsy = '60
  \font \eighttt = cmtt8 \hyphenchar\eighttt = -1
  \font \msbm = msbm10
  \font \sixbf = cmbx6
  \font \sixi = cmmi6 \skewchar \sixi = '177
  \font \sixrm = cmr6
  \font \sixsy = cmsy6 \skewchar \sixsy = '60
  \font \tensc = cmcsc10
  
  \font \titlefont = cmr7 scaled \magstep4
  \scriptfont \bffam = \bbseven
  \scriptscriptfont \bffam = \bbfive
  \textfont \bffam = \bbten

  \font\rs=rsfs10 

  \newskip \ttglue

  \def \eightpoint {\def \rm {\fam0 \eightrm }%
  \textfont0 = \eightrm
  \scriptfont0 = \sixrm \scriptscriptfont0 = \fiverm
  \textfont1 = \eighti
  \scriptfont1 = \sixi \scriptscriptfont1 = \fivei
  \textfont2 = \eightsy
  \scriptfont2 = \sixsy \scriptscriptfont2 = \fivesy
  \textfont3 = \tenex
  \scriptfont3 = \tenex \scriptscriptfont3 = \tenex
  \def \it {\fam \itfam \eightit }%
  \textfont \itfam = \eightit
  \def \sl {\fam \slfam \eightsl }%
  \textfont \slfam = \eightsl
  \def \bf {\fam \bffam \eightbf }%
  \textfont \bffam = \eightbf
  \scriptfont \bffam = \sixbf
  \scriptscriptfont \bffam = \fivebf
  \def \tt {\fam \ttfam \eighttt }%
  \textfont \ttfam = \eighttt
  \tt \ttglue = .5em plus.25em minus.15em
  \normalbaselineskip = 9pt
  \def \MF {{\manual opqr}\-{\manual stuq}}%
  \let \sc = \sixrm
  \let \big = \eightbig
  \setbox \strutbox = \hbox {\vrule height7pt depth2pt width0pt}%
  \normalbaselines \rm }


  \def \Headlines #1#2{\nopagenumbers
    \voffset = 2\baselineskip
    \advance \vsize by -\voffset
    \headline {\ifnum \pageno = 1 \hfil
    \else \ifodd \pageno \tensc \hfil \lcase {#1} \hfil \folio
    \else \tensc \folio \hfil \lcase {#2} \hfil
    \fi \fi }}

  \def \Title #1{\vbox{\baselineskip 20pt \titlefont \noindent #1}}

  \def \Date #1 {\footnote {}{\eightit Date: #1.}}

  \def \Authors #1{\bigskip \bigskip \noindent #1}

  \long \def \Addresses #1{\begingroup \eightpoint \parindent0pt
\medskip #1\par \par \endgroup }

  \long \def \Abstract #1{\begingroup \eightpoint
  \bigskip \bigskip \noindent
  {\sc ABSTRACT.} #1\par \par \endgroup }


  \def \lcase #1{\edef \auxvar {\lowercase {#1}}\auxvar }
  \def \vg #1{\ifx #1\null \null \else
    \ifx #1\ { }\else
    \ifx #1,,\else
    \ifx #1..\else
    \ifx #1;;\else
    \ifx #1::\else
    \ifx #1''\else
    \ifx #1--\else
    \ifx #1))\else
    { }#1\fi \fi \fi \fi \fi \fi \fi \fi \fi }

  \def \goodbreak {\vskip0pt plus.1\vsize \penalty -250 \vskip0pt
plus-.1\vsize }

  \newcount \secno \secno = 0
  \newcount \stno

  \def \seqnumbering {\global \advance \stno by 1
    \number \secno .\number \stno }

  \def \label #1{\def\localvariable {\number \secno
    \ifnum \number \stno = 0\else .\number \stno \fi }\global \edef
    #1{\localvariable }}

  \def\section #1{\global\def\SectionName{#1}\stno = 0 \global
\advance \secno by 1 \bigskip \bigskip \goodbreak \noindent {\bf
\number \secno .\enspace #1.}\medskip \noindent \ignorespaces}

  \long \def \sysstate #1#2#3{\medbreak \noindent {\bf \seqnumbering
.\enspace #1.\enspace }{#2#3\vskip 0pt}\medbreak }
  \def \state #1 #2\par {\sysstate {#1}{\sl }{#2}}
  \def \definition #1\par {\sysstate {Definition}{\rm }{#1}}

  \def \proof {\medbreak \noindent {\it Proof.\enspace }}
  \def \proofend {\ifmmode \eqno \square \else \hfill \square
\looseness = -1 \medbreak \fi }

  \def \$#1{#1 $$$$ #1}
  \def\=#1{\buildrel (#1) \over =}

  \def\Item #1{\smallskip \item {#1}}
  \newcount \zitemno \zitemno = 0
  \def\izitem {\zitemno = 0}
  \def\zitem {\global \advance \zitemno by 1 \Item {{\rm(\romannumeral
\zitemno)}}}

  \newcount \footno \footno = 1
  \newcount \halffootno \footno = 1
  \def\footcntr {\global \advance \footno by 1
  \halffootno =\footno
  \divide \halffootno by 2
  $^{\number\halffootno}$}

  \def \({\left (}
  \def \){\right )}
  \def \[{\left \Vert }
  \def \]{\right \Vert }
  \def \*{\otimes }
  \def \+{\oplus }
  \def \:{\colon }
  \def \<{\left \langle }
  \def \>{\right \rangle }
  \def \text #1{\hbox {\rm #1}}
  \def \and {\hbox {,\quad and \quad }}
  
  \def \calcat #1{\,{\vrule height8pt depth4pt}_{\,#1}}

  \def \crossproduct {{\hbox {\msbm o}}}
  \def \cstar {$C^*$}
  
  \def \inv {^{-1}}
  
  \def \square {\hbox {$\sqcap \!\!\!\!\sqcup $}}
  \def \stress #1{{\it #1}\/}
  
  \def \x {\times }
  \def \|{\Vert }
  \def \inv {^{-1}}


  \def\cite #1{{\rm [\bf #1\rm ]}}
  \def\scite #1#2{\cite{#1{\rm \hskip 0.7pt:\hskip 2pt #2}}}
  \def\lcite #1{(#1)}
  
  \def\bibitem#1#2#3#4{\smallskip \item {[#1]} #2, ``#3'', #4.}

  \def \references {
    \begingroup
    \bigskip \bigskip \goodbreak
    \eightpoint
    \centerline {\tensc References}
    \nobreak \medskip \frenchspacing }

  \font\titlefont = cmr7 scaled \magstep4
  \font\rs=rsfs10

  \def\({\left ( \vrule height 9pt width 0pt}
  \def\ATCP {{\curly T}}
  \def\B {{\curly B}}
  \def\F {{\bf F}}
  \def\Gen{{\cal G}}
  \def\I {{\cal I}}
  \def\J {{\cal J}}
  \def\Ker{{\rm Ker}}
  \def\Im{{\rm Im}}
  \def\Lin {{\curly L}}
  \def\Mat{{\cal A}}
  \def\N{{\bf N}}
  \def\Z{{\bf Z}}
  \def\OA{{\cal O}_\Mat}
  \def\O{{\cal O}}
  \def\R{{\cal R}}
  \def\T {{\curly T}}
  \def\TopExt{\ATCP(A,\a,\Tr)}
  \def\Tr{{\cal L}}
  \def\Univ {{\curly U}}
  \def\X{\Omega_\Mat}
  \def\a{\alpha}
  \def\bitem{\smallskip\item{$\bullet$}}
  \def\b{\beta}
  \def\cp{\mathop{\crossproduct_{\a,\Tr}} \N}
  \def\curly#1{\hbox{\rs #1\/}}
  \def\forro{,\ \ \forall}
  \def\half{^{1/2}}
  \def\map{\sigma}
  \def\scp{\crossproduct \N}

  \def\compos{\mathop{\raise 1pt \hbox{$\scriptscriptstyle \circ$}}}

  \newcount \bibno \bibno =0
  \def \newbib #1{\global \advance \bibno by 1 \edef #1{\number \bibno
}}

  \newbib \Baladi
  \newbib \Blackadar
  \newbib \Cuntz
  \newbib \CuntzTwo
  \newbib \CK
  \newbib \Deaconu
  \newbib \DeaconuMuhly
  \newbib \DR
  \newbib \amena
  \newbib \infinoa
  \newbib \AnHR
  \newbib \Murphy
  \newbib \Paschke
  \newbib \Pim
  \newbib \PV
  \newbib \QR
  \newbib \Renault
  \newbib \Ruelle
  \newbib \Stacey

  \def\titletext{A New Look at The Crossed-Product of a C*-algebra
\hfil\break by an Endomorphism}

  \Headlines {The Crossed-Product of a C*-algebra by an Endomorphism}
{Ruy Exel}

  \Title
  {\titletext}

  \Date {December 11, 2000}

  \Authors
  {Ruy Exel\footnote{*}{\eightrm Partially supported by CNPq.}}

  \Addresses
  {Departamento de Matem\'atica\par
  Universidade Federal de Santa Catarina\par
  88040-900 Florian\'opolis SC\par
  BRAZIL\vskip 4pt
  E-mail: exel@mtm.ufsc.br}

  \Abstract {We give a new definition for the crossed-product of a
\cstar-algebra $A$ by a *-endomorphism $\a$, which depends not only on
the pair $(A,\a)$ but also on the choice of a \stress{transfer
operator} (see definition below).  With this we generalize some of the
earlier constructions in the situations in which they behave best
(e.g.~for monomorphisms with hereditary range), but we get a different
and perhaps more natural outcome in other situations.  For example, we
show that the Cuntz--Krieger algebra $\OA$ arises as the result of our
construction when applied to the corresponding Markov subshift and a
very natural transfer operator.
  }

  \section{Introduction}
  Since the appearance of \cite{\Cuntz}, a paper of fundamental
importance in which Cuntz showed that $\O_n$ can be regarded as the
crossed-product of a UHF algebra by an endomorphism, many authors,
notably Paschke \cite{\Paschke}, Cuntz \cite{\CuntzTwo}, Doplicher and
Roberts \cite{\DR}, Stacey \cite{\Stacey}, and Murphy \cite{\Murphy},
have proposed general theories of crossed-products of \cstar-algebras
by single as well as semigroups of endomorphisms.  Most of these
theories are specially well adapted to deal with a monomorphism
(injective endomorphism) whose range is sometimes required to be a
hereditary subalgebra.

Motivated by Renault's Cuntz groupoid \scite{\Renault}{2.1}, Deaconu
introduced a construction \cite{\Deaconu} that can be applied to some
situations in which the endomorphism preserves the unit and therefore
usually does not have hereditary range.  This construction, and a
subsequent generalization in collaboration with Muhly
\cite{\DeaconuMuhly}, is designed to be applied to monomorphisms of
commutative \cstar-algebras, requiring special topological properties
of the map on the spectrum (branched covering) which are difficult to
phrase in noncommutative situations.

In this work we propose a new definition of crossed-product that
applies to most *-endomorphisms of unital \cstar-algebras,
generalizing some of the earlier constructions in the situations in
which they behave best (e.g.~monomorphisms with hereditary range), but
giving a different and perhaps more natural outcome in some other
situations.

The main novelty it that our construction depends not only on the
automorphism $\a$ of the \cstar-algebra $A$ but also on the choice a
\stress{transfer operator}, i.e.~a positive continuous linear map
$\Tr:A\to A$ such that $\Tr\(\a(a)b\) = a\Tr(b)$, for all $a,b\in A$.

Our motivation for considering transfer operators comes from classical
dynamics.  To be specific let $X$ be a compact topological space and
$\map:X\to X$ a continuous map, so that the pair $(X,\map)$ is a
classical dynamical system.  Suppose, in addition, that $\map$ is
surjective and a local homeomorphism.  Then the inverse image
$\map\inv(\{y\})$ is finite and nonempty for each $y\in X$ and the map
  $
  \Tr : C(X) \to C(X)
  $
  defined by
  $$
  \Tr (f)\calcat y = {1 \over \#\map\inv(\{y\})}
  \sum_{
  \buildrel {\scriptstyle t\in X} \over
  {\scriptstyle \map(t)=y}}
  f(t),\qquad \forall f\in C(X),\ y\in X,
  $$
  is often called a transfer operator for the pair $(X,\map)$ and has
important uses in classical dynamics (\cite{\Ruelle}, \cite{\Baladi}).

The realization that a transfer operator is a necessary ingredient for
the construction of the crossed-product comes from the fact that in
virtually all examples successfully treated so far the
crossed-product, say of the \cstar-algebra $A$ by the endomorphism
$\a$, turns out to be a \cstar-algebra $B$ which contains $A$ as well
as an extra element $S$ such that
  $$
  S^*AS\subseteq A.
  $$
  Since this takes place inside $A$, which is an ingredient of the
construction rather than the outcome, it is only natural that one
should be required to specify in advance what the map
  $$
  \Tr: a\in A \mapsto S^*aS \in A
  $$
  (which is necessarily positive) should be.  Again the examples seem
to warrant the requirement that
  $$
  Sa=\a(a)S\forro a\in A,
  $$
  and hence the axiom that $\Tr\(\a(a)b\) = a\Tr(b)$ must be accepted.
If we decide to demand no further properties of our extra element $S$
it makes sense to consider, at this point, the universal
\cstar-algebra generated by a copy of $A$ and an element $S$ subject
to the above restrictions. Thanks to Blackadar's theory of
\cstar-algebras given by generators and relations \cite{\Blackadar}
the universal algebra exists.  We call it $\TopExt$.

 In case $\a$ is an automorphism of $A$ and $\Tr$ is chosen to be the
inverse of $\a$ then the fact that $\Tr(1)=1$ implies that $S$ must be
an isometry and $\TopExt$ turns out to be the Toeplitz extension of
the usual crossed-product considered by Pimsner and Voiculescu in
\cite{\PV}.

  Since we are aiming at the crossed-product rather than the Toeplitz
extension, it is clear that the algebra that we seek is not yet
$\TopExt$, but is likely to be a quotient of it.  Quotients of
universal algebras are usually produced by throwing more relations at
it and hence we see that we must go in search of new relations.  If we
look at the literature for inspiration we will undoubtedly find
numerous occurrences of the relation
  $$
  \a(a) = S a S^*,
  \eqno{(\seqnumbering)}
  \label \BadEquation
  $$
  in which $S$ is moreover required to be an isometry,
  but this has some serious draw backs, among them that only an
endomorphism with hereditary range may be so implemented.  Stacey
\cite{\Stacey} suggests
  $$
  \a(a) = \sum_{i=1}^n S_i a S_i^*,
  $$
  but then one has to choose in advance which $n$ to pick.  Stacey
shows \scite{\Stacey}{2.1.b} that $n$ is not
intrinsic to the system $(A,\a)$ and hence there seems to be no 
good recipe for the choice of the number of terms in the above sum.

What seems to be the answer is a somewhat subtle point which, in
essence, boils down to a very interesting idea of Pimsner, namely the
passage from the Toeplitz--Cuntz--Pimsner algebra $\T_E$ to the
Cuntz--Pimsner algebra $\O_E$, in the context of Hilbert bimodules
(see \cite{\Pim} making sure to compare Theorems \lcite{3.4} and
\lcite{3.12}).

In order to describe Pimsner's idea, as adapted to the present
context, consider the subspace $M$ of $\TopExt$ obtained by taking the
closure of $AS$.  Using the defining properties of $S$ one easily
verifies that
  \bitem $MA \subseteq M$,
  \bitem $AM \subseteq M$, and
  \bitem $M^*M \subseteq A$.
  \medskip\noindent
  So $M$ may be viewed as a Hilbert bimodule over $A$ and, in
particular, $\overline{MM^*}$ is a \cstar-algebra.
  It follows that $M$ is invariant by the left multiplication
operators given by elements, either from $A$ or from
$\overline{MM^*}$.  It may therefore happen that for a certain $a\in
A$ and a certain $k\in\overline{MM^*}$ these operators
agree.  Equivalently
  $$
  abS = kbS \forro b\in A.
  $$
  In this case we call the pair $(a,k)$ a \stress{redundancy}.

Pimsner's idea \scite{\Pim}{Theorem 3.12} may be interpreted in the
present situation as the need to ``eliminate the redundancies'',
i.e.~to mod out the ideal generated by the differences $a-k$.
However, due to the fact that many of Pimsner's main hypothesis are
lacking here, namely the right module structure needs not be full
neither must the left action of $A$ on $M$ be isometric, we find that
eliminating all redundancies is a much too drastic approach.  Our
examples seem to indicate that the redundancies $(a,k)$ which should
be modded out are those for which $a$ lies in the closed two-sided
ideal generated by the range of $\a$.

When $\a(1)=1$ this ideal is the whole of $A$ and so our approach is
indeed to mod out all redundancies.  But, in order to make our theory
compatible with the successful theories for monomorphisms with
hereditary range, we seem to be restricted to modding out only the
redundancies indicated above.
  I must nevertheless admit that, besides the case in which $\a$
preserves the unit, I consider this as tentative, especially given the
relatively small number of interesting examples so far available.

  We thus define $A\cp$ to be the quotient of $\TopExt$ by the closed
two-sided ideal generated by the set of differences $a-k$, for all
redundancies $(a,k)$ such that $a\in \overline{A\R A}$, where $\R$
denotes the range of $\a$.

  It may be argued that our construction is not intrinsic to the
\cstar-dynamical system $(A,\a)$ since one needs to provide a transfer
operator, which is often not unique, by hand.  However
  in most of the examples that we have come across there seems to be a
transfer operator crying out to be picked.  For example, when dealing
with a monomorphism $\a$ of a \cstar-algebra $A$ with hereditary
range, the map
  $$
  E: a\in A \mapsto \a(1)\,a\,\a(1)\in \R,
  $$
  is a conditional expectation onto $\R$ and the composition $\Tr =
\a\inv \compos E$ is a transfer operator for $(A,\a)$.  Exploring this
example further we are able to show that, for any $a\in A$, the pair $\(\a(a),
SaS^*\)$ is a redundancy and hence \lcite{\BadEquation} does hold in
$A\cp$.  Moreover, since $\Tr(1)=1$, one has that $S$ an isometry
(which is not always the case).  We are then able to show that $A\cp$
is isomorphic to the universal \cstar-algebra for an isometry
satisfying equation \lcite{\BadEquation} and hence our theory is shown
to include Murphy's notion of crossed-product in the case of a single
monomorphism with hereditary range (see also \cite{\CuntzTwo}).
Likewise we show that Paschke's crossed-product algebras
\cite{\Paschke} may be recovered from our theory.

Leaving the realm of monomorphisms with hereditary range we next
discuss the case of a Markov (one-sided) subshift $(\X,\map)$ for a
finite transition matrix $\Mat$.  This clearly gives rise to the
endomorphism
  $$
  \a: f\in C(\X) \mapsto f\compos\map \in C(\X)
  $$
  which is unital and hence usually does not have hereditary range.

Here too one finds a very natural transfer operator and the
corresponding crossed-product is shown to be the Cuntz--Krieger
algebra $\O_A$.

  It is well known that $\map$ is surjective if and only if no column
of $\Mat$ is zero, which in turn is equivalent to $\a$ being
injective.
  Markov subshifts whose transition matrix possess trivial columns are
admittedly not among the most interesting ones but our result holds
just the same (with an Huef and Raeburn's definition of $\OA$
\cite{\AnHR}), thus showing that the present theory deals equally well
with endomorphisms which are not injective, as long as there is a
nontrivial transfer operator.

  The main problem left unresolved is the determination of the precise
conditions under which the natural map
  $$
  A \to A\cp
  \eqno{(\seqnumbering)}
  \label \IsItInjective
  $$
  is injective, although this is so in all of the examples considered.
We prove that the natural map
  $$
  A \to \TopExt
  $$
  is injective for any transfer operator whatsoever, but it will
certainly be necessary to assume special hypotheses in order for
\lcite{\IsItInjective} to be injective.  This is false, for instance,
if one chooses the zero map to play the role of the transfer operator.
But we conjecture that \lcite{\IsItInjective} is injective whenever
$\Tr$ is a non-degenerate transfer operator (see definition below).

  There are certainly numerous other examples where our theory may be
tested but in the sake of brevity we will leave these for a later
occasion.

  \section{Transfer operators}
  \label \RemarkTransfer
  Throughout this section, and most of this work, we will let $A$ be a
unital \cstar-algebra and
  $
  \a:A\to A
  $
  be a *-endomorphism.  The main concept to be introduced in this
section is as follows:

  \definition
  \label \TransferDefinitionOf
  A \stress{transfer operator} for the pair $(A,\a)$ is a continuous
linear map
  $$
  \Tr:A\to A
  $$
  such that
  \izitem
  \zitem $\Tr$ is positive in the sense that $\Tr(A_+)\subseteq A_+$,
and
  \zitem $\Tr\(\a(a)b\) = a\Tr(b)$, for all $a,b\in A$.

Let $\Tr$ be a transfer operator for $(A,\a)$.  As a consequence of
positivity it is clear that $\Tr$ is self-adjoint and hence the
simmetrized version of \lcite{\TransferDefinitionOf.ii}, namely
  $\Tr\(a\a(b)\) = \Tr(a)b$
  holds as well.  In particular this implies that the range of $\Tr$
is a two-sided ideal in $A$.  Also
  $$
  a\Tr(1) = \Tr(\a(a)) = \Tr(1)a,
  \eqno{(\seqnumbering)}
  \label \TransferAtOne
  $$
  for all $a\in A$,
  so $\Tr(1)$ is a positive central element in $A$.  In our examples
we will see that quite often, though not always, one has $\Tr(1)=1$,
in which case $\Tr$ will necessarily be surjective.

It is interesting to notice that if $a\in\Ker(\a)$ (observe that we
are not assuming $\a$ to be injective), and $b\in A$ we have
  $$
  a\Tr(b) = \Tr(\a(a)b) = 0,
  $$
  so $\Ker(\a)$ is orthogonal to the range of $\Tr$.  Hence, the less
injective $\a$ is, the less room there is for a nontrivial transfer
operator to exist (notice that the zero map is always a transfer
operator).

Given the predominant role to be played by the range of $\a$ we will
adopt the notation
  $$
  \R := \a(A).
  $$
  It is well known that $\R$ is a sub-\cstar-algebra of $A$.

  Define the map $E:A \to A$ by $E = \a\compos\Tr$,
  and observe that if $b\in\R$, say $b = \a(c)$ with $c\in A$, one has
for all $a\in A$ that
  $$
  E(ab) =
  \a(\Tr(a\a(c))) =
  \a(\Tr(a)c) =
  E(a)\a(c) =
  E(a)b,
  $$
  and similarly $E(ba) = bE(a)$.  Therefore $E$ is linear for the
natural $\R$--bimodule structure of $A$.

Clearly $E$ maps $A$ into $\R$ so one sees that $E$ satisfies almost
all of the axioms of a conditional expectation.

  \state Proposition
  \label \DefineNonDegenerateTransfer
  Let $\Tr$ be a transfer operator for the pair $(A,\a)$.  Then the
following are equivalent:
  \izitem
  \zitem the composition $E = \a\compos\Tr$ is a conditional
expectation onto $\R$,
  \zitem $\a\compos\Tr\compos\a =\a$,
  \zitem $\a(\Tr(1))=\a(1).$

  \proof
  Assuming (i) one has that $E|_\R = id_\R$ and hence for all $a\in
A$,
  $$
  \a(\Tr(\a(a))) = E(\a(a)) = \a(a),
  $$
  proving (ii).
  Observe that plugging $a=1$ in \lcite{\TransferAtOne} we get
  $
  \Tr(1) = \Tr(\a(1)).
  $
  Thus, if (ii) is assumed we have
  $$
  \a(1) = \a(\Tr(\a(1))) = \a(\Tr(1)).
  $$
  Supposing that (iii) holds we have for all $a\in A$ that
  $$
  \a(\Tr(\a(a))) =
  \a(\Tr(\ \a(a)1\ )) =
  \a(a\Tr(1)) =
  \a(a)\a(\Tr(1)) =
  \a(a)\a(1) =
  \a(a),
  $$
  and (ii) follows.  Finally it is easy to see that (ii) implies that
$E$ is the identity on $\R$.  From this and from our discussion above
we deduce (i).
  \proofend

  This motivates our next:

  \definition
  We will say that a transfer operator $\Tr$ is
\stress{non-degenerate} if the equivalent conditions of
\lcite{\DefineNonDegenerateTransfer} hold.

  \state Proposition
  Let $\Tr$ be a non-degenerate transfer operator.  Then $A$ may be
written as the direct sum of ideals
  $$
  A = \Ker(\a) \+ \Im(\Tr).
  $$
  
  \proof
  Regardless of $\Tr$ being non-degenerate we have seen that
$\Ker(\a)$ is orthogonal to the range of $\Tr$.  For all $a\in A$ we
have that
  $$
  a = \(a-\Tr\a(a)\) + \Tr\a(a)
  $$
  and one may easily check that if $\Tr$ is non-degenerate then
$a-\Tr\a(a)$ lies in $\Ker(\a)$.
  \proofend

Let us now discuss a fairly general way to produce non-degenerate
transfer operators.

  \state Proposition
  \label \ConstrNDTransfer
  Suppose that $E$ is a positive conditional expectation from $A$ onto
$\R$.
  Suppose also that we are given an ideal $\J$ in $A$ such that $A =
\J \+ \Ker(\a)$ (e.g.~$\J=A$ when $\a$ is injective).  Let
$\b:\R\to\J$ be the inverse of the restriction $\a|_{\J} : \J \to \R$.
Then the composition
  $
  \Tr = \b \compos E
  $
  is a non-degenerate transfer operator for $(A,\a)$ which moreover
satisfies
  $\a\compos\Tr=E$.

  \proof
  It is clear that $\Tr$ is positive.  For $a,b\in A$ we have
  $$
  \Tr(\a(a)b) =
  \b(E(\a(a)b)) =
  \b(\a(a)E(b)).
  $$
  We next claim that the last term above coincides with $a\b(E(b))$.
In order to see this observe that these are both in $\J$ where the
restriction of $\a$ is 1--1.  Hence the equality sought is equivalent
to
  $$
  \a\Big(\ \b(\a(a)E(b))\ \Big) =
  \a\Big(\ a\b(E(b))\ \Big),
  $$
  which the reader may verify with no difficulty using that
  $\a\compos\b$ is the identity on $\R$.
  Finally observe that
  $$
  \a\compos\Tr =
  \a\compos\b\compos E =
  E,
  $$
  proving that $\Tr$ satisfies \lcite{\DefineNonDegenerateTransfer.i}
and hence is non-degenerate.
  \proofend

  \section{The Crossed-Product}
  \label \ToeplitzSection
  As before we will let $A$ be a \cstar-algebra and $\a:A\to A$ be a
*-endomorphism.  From now on we will also fix a transfer operator
$\Tr$ for the pair $(A,\a)$.
  Even though for a degenerate $\Tr$ (e.g.~$\Tr=0$) the
crossed-product is likely not to be very interesting, it seems that
our formalism may be carried out to a certain extent without any
further hypothesis on $\Tr$.
  In the examples we will explore, however, $\Tr$ will often be
obtained by an application of \lcite{\ConstrNDTransfer} and hence will
be non-degenerate.

Our next immediate goal is to introduce a \cstar-algebra which will be
an extension of the crossed-product algebra to be defined later.

  \definition
  \label \DefineToeplitzExtension
  Given a *-endomorphism $\a$ of a unital \cstar-algebra $A$ and a
transfer operator $\Tr$ for $(A,\a)$ we let
  $\TopExt$
  be the universal unital \cstar-algebra generated by a copy of $A$
and an element $S$ subject to the relations
  \izitem
  \zitem $Sa = \a(a)S$,
  \zitem $S^*aS = \Tr(a)$,
  \medskip\noindent
  for every $a\in A$.

For any representation of the above relations one obviously has that
$\|S\| = \|\Tr(1)\|^{1/2}$, so these are \stress{admissible relations}
in the sense of \cite{\Blackadar}, from where one deduces the
existence of $\TopExt$.

When $\Tr(1)=1$ it follows that $S$ is an isometry but, contrary to
some existing notions of crossed-products (\cite{\Murphy},
\cite{\Paschke}), it might well be that $S$ is not even a partial
isometry.  Another difference with respect to the existing literature
is that we do not assume that $SaS^*=\a(a)$, or even that
$SAS^*\subseteq A$.  Such conditions would knock out some of our most
interesting examples.

An essential feature of the description of $\TopExt$ as a universal
\cstar-algebra is a canonical map
  $$
  A\to \TopExt
  $$
  (see \cite{\Blackadar}),
  but there is no a priori guarantee that this will be injective.  As
is the case with most objects defined by means of generators and
relations there may or may not be hidden relations which force nonzero
elements of $A$ to be zero in the universal object.  It is our next
main goal to prove that in the present situation $A$ embeds in
$\TopExt$ injectively.

Let $A_\Tr$ be a copy of the underlying vector space of $A$.  Define
the structure of a right $A$--module on $A_\Tr$ by
  $$
  m\cdot a = m\a(a) \forro m\in A_\Tr,\ a\in A.
  $$
  Also define an $A$-valued (possibly degenerate) inner-product
  $
  \<\cdot\,,\cdot\>_\Tr
  $
  on $A_\Tr$ by
  $$
  \<m,n\>_\Tr=\Tr(m^*n)\forro m,n\in A_\Tr.
  $$
  Upon modding out vectors of norm zero and completing we get a right
Hilbert $A$--module which we denote by $M_\Tr$.
  For every $a$ in $A$ and $m$ in $A_\Tr$ we have
  $$
  \|\<am,am\>\| =
  \|\Tr(m^*a^*am)\| \leq
  \|a\|^2\|\Tr(m^*m)\| =
  \|a\|^2\|\<m,m\>\|,
  $$
  so we see that left multiplication by $a$ on $A_\Tr$ extends to a
bounded operator
  $
  \phi(a) : M_\Tr \to M_\Tr,
  $
  which may be easily proven to be adjointable.  This in turn defines
a *-homomorphism
  $$
  \phi:A \to \Lin(M_\Tr).
  $$
  thus making $M_\Tr$ a Hilbert bimodule over $A$ (in the sense of
\cite{\Pim}).

  \state Lemma
  If $\Tr$ is a transfer operator for $(A,\a)$ and $n\in \N\cup\{0\}$
then $\Tr^n$ is a transfer operator for $(A,\a^n)$.  Moreover the map
  $$
  \gamma_n:
  x\in A_{\Tr^n} \ \longmapsto\ \a(x) \in A_{\Tr^{n+1}}
  $$
  extends to an adjointable linear map
  $$
  \gamma_n: M_{\Tr^n} \to M_{\Tr^{n+1}},
  $$
  with $\|\gamma_n\|\leq \|\Tr(1)\|\half$, whose adjoint is given by
  $$
  \gamma_n^*(x)=\Tr(x)
  \forro x\in A_{\Tr^{n+1}}.
  $$

  \proof
  Clearly $\Tr^n$ is positive.  With respect to
\lcite{\TransferDefinitionOf.ii} let $a,b\in A$.  By induction we have
  $$
  \Tr^n\(\a^n(a)b\) =
  \Tr^{n-1}\Big(\Tr\(\a(\a^{n-1}(a))b\)\Big) =
  \Tr^{n-1}\( \a^{n-1}(a) \Tr(b)\) =
  a \Tr^{n-1}\( \Tr(b)\) =
  a \Tr^n(b),
  $$
  proving that $\Tr^n$ is in fact a transfer operator for $(A,\a^n)$.
Next observe that for all $x\in A_{\Tr^{n}}$ one has
  $$
  \<\a(x),\a(x)\>_{\Tr^{n+1}} =
  \Tr^n\(\Tr(\a(x^*)\a(x))\) =
  \Tr^n(x^*\Tr(1)x) \$\leq
  \|\Tr(1)\|\ \Tr^n(x^*x) =
  \|\Tr(1)\|\ \<x,x\>_{\Tr^n}.
  $$
  It follows that $\gamma_n$ is bounded on $A_{\Tr^{n}}$ and hence
extends to $M_{\Tr^{n}}$ with $\|\gamma_n\|\leq \|\Tr(1)\|\half$.

  Let $x\in A_{\Tr^{n}}$ and $y\in A_{\Tr^{n+1}}$.  Then
  $$
  \<\gamma_n(x),y\>_{\Tr^{n+1}} =
  \Tr^n\(\Tr(\a(x^*)y)\) =
  \Tr^n\(x^*\Tr(y)\) =
  \<x,\Tr(y)\>_{\Tr^n}.
  $$
  We therefore have that
  $$
  \|\Tr(y)\|_{\Tr^n}^2 =
  \|\<\Tr(y),\Tr(y)\>_{\Tr^n}\| =
  \|\<\gamma_n(\Tr(y)),y\>_{\Tr^{n+1}}\| \leq
  \|\gamma_n\|\,\|\Tr(y)\|_{\Tr^n}\|y\|_{\Tr^{n+1}},
  $$
  which implies that
  $\|\Tr(y)\|_{\Tr^n}\leq \|\gamma_n\|\,\|y\|_{\Tr^{n+1}}$,
  and hence the correspondence $y\to\Tr(y)$ extends continuously to a
map from $A_{\Tr^{n+1}}$ to $A_{\Tr^{n}}$ which may now easily be
shown to be the adjoint of $\gamma_n$.
  \proofend

  The following should now be evident:

  \state Proposition
  \label \MInfinite
  Consider the Hilbert bimodule $M_\infty$ over $A$ given by
  $$
  M_\infty = \bigoplus_{n=0}^\infty M_{\Tr^n},
  $$
  and define
  $
  S : M_\infty \to M_\infty
  $
  by
  $$
  S(x_0,x_1,\ldots) = \(0,\gamma_0(x_0),\gamma_1(x_1),\ldots\).
  $$
  Then
  $S$ is an adjointable map on $M_\infty$ and its adjoint is given by
  $$
  S^*(x_0,x_1,x_2,\ldots) = \(\Tr(x_1),\Tr(x_2),\ldots\),
  $$
  whenever $x = (x_0,x_1,\ldots) \in M_\infty$ is such that $x_n\in
A_{\Tr^n}$ for every $n$.

  \proof
  Left to the reader.
  \proofend

  We now come to a main technical result:

  \state Theorem
  Let $\a$ be a *-endomorphism of a unital \cstar-algebra $A$ and let
$\Tr$ be a transfer operator for $(A,\a)$.  Then there exists a
faithful representation $\rho$ of $A$ on a Hilbert space $H$ and an
operator $S\in\B(H)$ such that
  \izitem
  \zitem $S\rho(a) = \rho(\a(a))S$,
  \zitem $S^*\rho(a)S = \rho(\Tr(a))$,
  \medskip\noindent
  for every $a\in A$.

  \proof
  Consider the \cstar-algebra $\Lin(M_\infty)$ and view the left
$A$--module structure of $M_\infty$ as *-homomor\-phism
  $
  \rho:A \to \Lin(M_\infty).
  $
  Observe that each $M_{\Tr^n}$ is invariant under $\rho$ and that the
restriction of $\rho$ to $M_{\Tr^0}$ is faithful.  It follows that
$\rho$ itself is faithful.
  Considering the operator $S$ defined in \lcite{\MInfinite} one may
now easily prove that
  $$
  S\rho(a) = \rho(\a(a))S \and S^*\rho(a)S = \rho(\Tr(a)) \forro a\in
A.
  $$
  It now suffices to compose $\rho$ with any faithful Hilbert space
representation of $\Lin(M_\infty)$.
  \proofend

  As an immediate consequence we have:

  \state Corollary
  \label \ResumeFiel
  Let $\a$ be a *-endomorphism of a unital \cstar-algebra $A$ and let
$\Tr$ be a transfer operator for $(A,\a)$.  Then
  the canonical map
  $
  A\to \TopExt
  $
  is injective.

From now on we will view $A$ as a sub-\cstar-algebra of $\TopExt$, by
\lcite{\ResumeFiel}.  The canonical representation of $S$ within
$\TopExt$ will be denoted simply by $S$.

  \definition
  (Compare \scite{\Pim}{Theorem 3.12}).
  By a \stress{redundancy} we will mean a pair $(a,k)\in A\times
\overline{ASS^*A}$ such that
  $$
  abS = kbS,
  $$
  for all $b\in A$.

We may now introduce our main concept.  In what follows we will be
mostly interested in the redundancies $(a,k)$ for $a\in \overline{A\R
A}$ (the closed two-sided ideal generated by $\R = \a(A)$).

  \definition
  \label \CrossedProductDef
  The crossed-product of the unital \cstar-algebra $A$ by the
*-endomorphism $\a:A\to A$ relative to the choice of the transfer
operator $\Tr$, which we denote by $A\cp$, or simply by $A\scp$ when
$\a$ and $\Tr$ understood, is defined to be the quotient of $\TopExt$
by the closed two-sided ideal $\I$ generated by the set of differences
$a-k$, for all redundancies $(a,k)$ such that $a\in \overline{A\R A}$.

Unfortunately we have not been able to show that the natural inclusion
of $A$ in $A\cp$ is injective.  We will see, however, that this is the
case in a series of examples which we will now discuss.

  \section{Monomorphisms with hereditary range}
  We now wish to study the crossed-product under special hypotheses
which have been considered in \cite{\Murphy}.
  Let $A$ be a unital \cstar-algebra and let $\a:A\to A$ be a
*-endomorphism.  Denote by $P$ the idempotent $\a(1)$ and notice that
the range of $\a$, which we keep denoting by $\R$, is contained in the
hereditary subalgebra $PAP$.

  \state Proposition
  \label \HereditaryIsPAP
  If $\R$ is a hereditary subalgebra of $A$ then $\R=PAP$.

  \proof
  As noticed above $\R \subseteq PAP$.  Given $a\in PAP$ with $a\geq
0$ we have that
  $$
  0 \leq a = PaP \leq \|a\| P \in \R,
  $$
  which implies that $a\in \R$.
  \proofend

  Throughout this section we will suppose that $\R$ is hereditary and
hence that $\R=PAP$.  It follows that the map
  $$
  E: a\in A \mapsto PaP \in \R
  $$
  is a conditional expectation onto $\R$.  We will assume moreover
that $\a$ is injective so that we are under the hypothesis of
\lcite{\ConstrNDTransfer} (with $\J=A$) and hence the composition
  $$
  \Tr = \a\inv\compos E
  \eqno{(\seqnumbering)}
  \label \MyTransfer
  $$
  defines a non-degenerate transfer operator for $(A,\a)$.
 
  \state Proposition
  \label \SIsAnIsometry
  Suppose that $\a$ is injective, that $\R$ is hereditary, and that
$\Tr$ is given as above.  Then:
  \izitem
  \zitem The canonical element $S\in \TopExt$ is an isometry, and
hence also its image $\dot S\in A\scp$.
  \zitem For every $a\in A$ one has that
  $
  \(\a(a),S a S^*\)
  $
  is a redundancy.

  \proof
  With respect to (i) we have by \lcite{\DefineToeplitzExtension.ii}
that
  $$
  S^*S = \Tr(1) = \a\inv(E(1)) = \a\inv(P) =\a\inv(\a(1)) = 1.
  $$
  In order to prove (ii) let $b\in A$ and notice that
  $$
  SaS^*bS =
  Sa\a\inv(E(b)) =
  \a(a)E(b)S =
  \a(a)PbPS =
  \a(a)bS,
  $$
  where we have used that $PS = \a(1) S = S 1 = S$.
  It remains to notice that
  $SaS^* = \a(a)SS^*\in \overline {A S S^* A}$.
  \proofend

  \definition
  We will denote by $\Univ(A,\a)$ the universal unital \cstar-algebra
generated by $A$ and an isometry $T$ subject to the relation $\a(a) =
TaT^*$, for every $a\in A$.

$\Univ(A,\a)$ has been proposed (\cite{\CuntzTwo}, \cite{\Stacey},
\cite{\Murphy}) as the definition for the crossed-product of $A$ by
$\a$.  It is easy to see that, under the present hypothesis that $\a$
is a monomorphism, $A$ embeds in $\Univ(A,\a)$ injectively.  In fact,
using the terminology of the last paragraph of the section entitled
``Examples of Simple \cstar-algebras'' in \cite{\CuntzTwo} (where $\a$
is denoted by $\varphi$), consider a faithful representation of
$A^\infty \crossproduct_{\a^\infty}\Z$ on a Hilbert space $H$, and let
$U$ be the unitary implementing the automorphism $\a^\infty$.  Then,
denoting the unit of $A$ by $1_A$, one may easily prove that $T:=
\a(1_A)U1_A$ is an isometry on $1_AH$, and that $\a(a) = TaT^*$, for
all $a$ in $A$.  By universality there exists a canonical
representation of $\Univ(A,\a)$ on $\Lin(1_AH)$ which is the identity
on $A$ and hence the canonical map $A \to \Univ(A,\a)$ must be
injective (see also \scite{\Stacey}{Section 2} and
\scite{\Murphy}{Section 2}).

We will therefore view $A$ as a subalgebra of $\Univ(A,\a)$.  As a
consequence of \lcite{\SIsAnIsometry} we have:
  
  \state Corollary
  \label \DeMurParaCP
  Under the hypotheses of \lcite{\SIsAnIsometry}
  there exists a unique *-epimomorphism
  $$
  \phi: \Univ(A,\a) \to A\cp
  $$
  such that
  $\phi(T)=\dot S$, and
  $\phi(a)=\dot a$, for all $a\in A$, where
  $\dot S$ and $\dot a$ are the canonical images of $S$ and $a$ in
$A\scp$, respectively.

  \proof
  By \lcite{\SIsAnIsometry.i} $\dot S$ is an isometry and by
\lcite{\SIsAnIsometry.ii} we have that $\a(\dot a) = \dot S \dot a
\dot S^*$.  Hence the conclusion follows from the universal property
of $\Univ(A,\a)$.
  \proofend

  On the other hand we have:

  \state Proposition
  \label \HereditaryCase
  Under the hypotheses of \lcite{\SIsAnIsometry}
  there exists a unique *-epimomorphism
  $$
  \psi: \TopExt \to \Univ(A,\a)
  $$
  such that
  $\psi(S) = T$,
  and $\psi(a)=a$, for all $a\in A$.

  \proof Observe that, for all $a\in A$,
  $$
  T a = T a T^*T = \a(a)T,
  $$
  and
  $$
  T^*aT =
  T^*TT^*aTT^*T =
  T^*PaPT =
  T^*E(a)T =
  T^*\a(\a\inv(E(a)))T \$=
  T^*T \Tr(a)T^*T =
  \Tr(a).
  $$
  In other words, relations \lcite{\DefineToeplitzExtension.i--ii}
hold for $A$ and $T$ within $\Univ(A,\a)$ and hence the conclusion
follows from the universal property of $\TopExt$.
  \proofend

  Putting together the conclusions of \lcite{\DeMurParaCP} and
\lcite{\HereditaryCase} we see that the quotient map from $\TopExt$ to
$A\cp$ factors through $\Univ(A,\a)$ yielding the commutative diagram
  \bigskip
  $$
  \TopExt\quad \longrightarrow\quad A\cp
  $$$$
  _{\displaystyle \psi} \searrow \kern 40pt \nearrow_{\displaystyle
\phi}
  $$$$
 \Univ(A,\a)
  $$
  \bigskip\noindent
  in which the horizontal arrow is the quotient map.
  In the main result of this section we will show that $\phi$ is in
fact an isomorphism.

  \state Theorem
  \label \MurphyAsCP
  Let $A$ be a unital \cstar-algebra and let $\a:A\to A$ be an
injective *-endomorphism whose range is a hereditary subalgebra of
$A$.  Consider the transfer operator $\Tr$ given by
  $$
  \Tr(a) = \a\inv(P a P)\forro a\in A,
  $$
  where $P=\a(1)$.  Then the map $\phi$ of \lcite{\DeMurParaCP} is a
*-isomorphism between $\Univ(A,\a)$ and $A\cp$.

  \proof
  We begin by claiming that the map $\psi$ of \lcite{\HereditaryCase}
vanishes on the ideal $\I$ mentioned in \lcite{\CrossedProductDef}.
In order to prove this let $(a,k)\in A\x \overline{ASS^*A}$ be a
redundancy with $a\in\overline{A\R A}$.  Therefore for all $b\in A$
one has
  $
  abS = kbS.
  $
  Applying $\psi$ to both sides of this equation gives
  $$
  a b T = \psi(k) b T.
  $$
  Observing that $P =\a(1) = TT^*$ we have for all $b,c\in A$ that
  $$
  a b P c =
  a b T T^* c =
  \psi(k) b T T^* c =
  \psi(k) b P c.
  $$
  It follows that $a x = \psi(k)x$ for all
  $ x\in \overline {APA}$.
  Since $k\in \overline{A SS^*A}$ we have that
  $\psi(k)\in \overline{A TT^*A} = \overline{APA}$.
  Finally we have that
  $$
  a \in
  \overline {A\R A} =
  \overline {APAPA} =
  \overline {APA}.
  $$
  Therefore we must have that $a=\psi(k)$.  It follows that
$\psi(a-k)=0$ and hence that $\psi$ vanishes on $\I$ as claimed.  By
passage to the quotient we get a map
  $$
  \widetilde\psi : A\scp \to \Univ(A,\a),
  $$
  which one may now easily prove to be the inverse of the map $\phi$
of \lcite{\DeMurParaCP}.
  \proofend

  \section{Paschke's Crossed-Product}
  Let us now comment on the relationship between Paschke's notion of
crossed-product \cite{\Paschke} and the one introduced in this work.
  Following \cite{\Paschke} we will fix, throughout this section, a
unital \cstar-algebra $A$ acting on a Hilbert space $H$, and an
isometry $T$ such that both $TAT^*$ and $T^*AT$ are contained in $A$.
Following Paschke's notation we let $C^*(A,T)$ be the \cstar-algebra
of operators on $H$ generated by $A$ and $T$.
  One may then define an endomorphism
  $$
  \a: A \to A
  $$
  by the formula $\a(a) = TaT^*$.
  Because $T^*\a(a)T = a$ we have that $\a$ is necessarily injective.
Let $$P:=\a(1)=TT^*,$$ and observe that the range of $\a$ is precisely
the hereditary subalgebra
  $
  \R = PAP
  $
  of $A$.
  This suggests choosing the conditional expectation
  $$
  E: a\in A \mapsto PaP\in \R,
  $$
  and the transfer operator
  $\Tr = \a\inv \compos E$, as in \lcite{\MyTransfer}.

  \state Theorem
  Let $A$ and $T$ be as above and suppose that the hypothesis of\/
\scite{\Paschke}{Theorem 1} are satisfied (that is, $A$ is strongly
amenable, $T$ is not unitary, and there is no nontrivial ideal $J$ of
$A$ such that $TJT^*\subseteq J$).  Then $C^*(A,T)$ is canonically
isomorphic to $A\cp$.

  \proof
  By \lcite{\MurphyAsCP} it follows that there exists a
*-epimomorphism
  $$
  \phi: A\cp \to C^*(A,T),
  $$
  such that $\phi(\dot a)=a$, for all $a\in A$, and $\phi(\dot S)=T$.
  We will now use \scite{\Paschke}{Theorem 1} in order to prove that
$A\scp$ is simple.  By \lcite{\SIsAnIsometry} we have that $\dot S$ is
an isometry and that
  $
  \dot S A \dot S^* = \a(A) \subseteq A,
  $
  while the fact that
  $
  \dot S^* A \dot S = \Tr(A) \subseteq A
  $
  follows from \lcite{\DefineToeplitzExtension.ii}.
  
  Given that $T$ is not unitary and that $\phi(\dot S)=T$, it is clear
that $\dot S$ is not unitary either.
  %
  %
  For any ideal $J$ of $A$ we have that $\dot S J\dot S^* = \a(J) = T
J T^*$, and hence there is no nontrivial ideal $J$ of $A$ such that
$\dot S J\dot S^*\subseteq J$.  By Paschke's Theorem
\scite{\Paschke}{Theorem 1} we have that $C^*(A,\dot S) = A\scp$ is
simple and hence that $\phi$ is an isomorphism.
  \proofend

  \section{Cuntz--Krieger algebras as Crossed-Products}
  Throughout this section we will let $\Mat$ be an $n\times n$ matrix
with $\Mat(i,j)\in\{0,1\}$ for all $i$ and $j$, and such that no row
of $\Mat$ is identically zero.
  Our main goal is to show that the Cuntz--Krieger algebra $\OA$ (see
\cite{\CK}, \cite{\AnHR}) may be characterized as the crossed-product
of a commutative \cstar-algebra by an endomorphism arising from Markov
sub-shifts.

For convenience we let
  $$
  \Gen = \{1,2,\ldots,n\}.
  $$
  Considering the Cantor space $\Gen^\N$ (we let $\N =
\{1,2,3,\ldots\}$ by convention) let $\X$ be the compact subspace of
$\Gen^\N$ given by
  $$
  \X = \{\xi=(\xi_i)_{i\in\N}\in\Gen^\N: \Mat(\xi_i,\xi_{i+1})=1,\
\forall i\in\N\}.
  $$
  Let $\map:\X\to \X$ be given by
  $$
  \map(\xi_1,\xi_2,\ldots) = (\xi_2,\xi_3,\ldots),
  $$
  so that $(\X,\map)$ is the well known dynamical system usually
referred to as the Markov sub-shift.
  Given $\xi\in\X$ observe that
  $$
  \map\inv(\{\xi\}) = \big\{(x,\xi_1,\xi_2,\ldots):x\in\Gen,\
\Mat(x,\xi_1)=1\big\}.
  $$
  In particular the number of preimages of $\xi$ is exactly the number
of elements $x\in\Gen$ such that $\Mat(x,\xi_1)=1$.  In other words
  $$
  \#\map\inv(\{\xi\}) = \sum_{x\in\Gen}\Mat(x,\xi_1).
  $$
  From now on we will let $Q$ be the integer-valued function on $\X$
defined by
  $$
  Q(\xi) = \sum_{x\in\Gen}\Mat(x,\xi_1).
  $$
  Since $Q$ depends only on the first component $\xi_1$ of $\xi$ we
have that $Q$ is a continuous function on $\X$.
  Observe also that $\map(\X) = \{\xi\in\X: Q(\xi) > 0\}$, whence
$\map(\X)$ is a clopen set.  In particular $\map$ is surjective if and
only if $Q>0$, which in turn is equivalent to the fact that no column
of $\Mat$ is zero.

  Consider the endomorphism of the commutative \cstar-algebra $C(\X)$
given by
  $$
  \a(f) = f\compos\map\forro f\in C(\X).
  $$
  Clearly $\a$ is injective if and only if $\map$ is surjective, which
we have seen to be equivalent to the absence of trivial columns.
  Markov sub-shifts with a transition matrix $\Mat$ possessing
identically zero columns are admittedly not among the most interesting
examples.  Nevertheless we will allow zero columns here mostly to
illustrate that our theory deals equally well with endomorphisms which
are not injective.

  The kernel of $\a$ therefore consists of the ideal
$C_0\(\X\setminus\map(\X)\)$ formed by the functions $f$ vanishing on
$\map(\X)$.  Since the latter is a clopen set we have that $C(\X)$ may
be written as the direct sum of ideals
  $$
  C_0\(\map(\X)\)\+\Ker(\a),
  $$
  as required by \lcite{\ConstrNDTransfer}.
  Accordingly we will let $\J = C_0\(\map(\X)\)$ and we will denote by
$\b$ the inverse of the map $\a|_{\J} : \J \to \R$, where $\R =
\a(C(\X))$, as usual.

  Given $f\in C(\X)$ let $E(f)$ be the function on $\X$ defined by
  $$
  E(f)\calcat {\xi} = {1 \over \#\map\inv(\{\map(\xi)\})}
  \sum_{
  \buildrel {\scriptstyle \eta\in \X} \over
  {\scriptstyle \map(\eta)=\map(\xi)}}
  f(\eta)\forro \xi\in \X,
  $$
  so that $E(f)\calcat \xi$ is the average of the values of $f$ on the
elements of $\X$ which have the same image as $\xi$ under $\map$.  We
leave it for the reader to prove:

  \state Proposition
  For every $f\in C(\X)$ one has that $E(f)\in C(\X)$.  Moreover
  $E$ is a positive conditional expectation from $C(\X)$ to $\R$.

  Following \lcite{\ConstrNDTransfer} we therefore have that
  $
  \Tr=\b\compos E
  $
  is a transfer operator for $(C(\X),\a)$.  One may alternatively
define $\Tr$ directly as
  $$
  \Tr(f)\calcat \xi =
  \left\{ \matrix{
    \displaystyle {1 \over
                  \#\map\inv(\{\xi\})}
    \sum_{
    \buildrel {\scriptstyle \eta\in \X} \over
    {\scriptstyle \map(\eta)=\xi}}
    f (\eta),
  & \hbox{if $\xi\in\map(\X)$,}
  \cr\cr\cr
  0, & \hbox{ if $\xi\notin\map(\X)$.}
  } \right.
  $$

  It is the goal of this section to prove:
  
  \state Theorem
  \label \CKIsomorphism
  For every $n\times n$ matrix $\Mat$ with no zero rows one has that
the Cuntz--Krieger \cstar-algebra $\OA$ is isomorphic to $C(\X)\cp$.

  \proof The proof will be accomplished in several steps.  We begin by
introducing some notation.  Given $x\in\Gen$ we let
  $Q_x$ and $P_x$ be the continuous functions defined for all $\xi\in
\X$ by
  $$
  Q_x(\xi) = \Mat(x,\xi_1)
  \and
  P_x(\xi) = [\xi_1=x],
  $$
  where the brackets correspond to the obvious boolean valued
function.
  Observe that the function $Q$ defined above may then be written as
  $
  Q = \sum_{x\in\Gen} Q_x.
  $
  Observe also that for all $x\in\Gen$ one has that
  $$
  Q_x = \sum_{y\in\Gen} \Mat(x,y)P_y.
  $$
  Working within $\ATCP(C(\X),\a,\Tr)$ let, for each $x\in\Gen$,
  $$
  S_x = P_x\a(Q)\half S,
  $$
  and observe that
  $$
  S_x^*S_x =
  S^*\a(Q) P_x S =
  \Tr\(\a(Q) P_x\) =
  Q_x,
  $$
  where the last step should be verified by direct computation.  Since
$Q_x$ is a projection we have that $S_x$ is a partial isometry.
  We now claim that for all $x\in\Gen$,
  $$
  S_x^*S_x \equiv \sum_{y\in\Gen} \Mat(x,y) S_yS_y^*
  \eqno{(\dagger)}
  $$
  modulo the ideal $\I$ of \lcite{\CrossedProductDef}.
  In order to see this observe that, for $y\in\Gen$ and all $b\in
C(\X)$, one has that
  $$
  S_y S_y^*bS =
  P_y\a(Q)\half S S^* \a(Q)\half P_y bS=
  P_y\a(Q)\half S \Tr\(\a(Q)\half P_y b\) \$=
  P_y\a(Q)\half E\(\a(Q)\half P_y b\) S=
  P_y\a(Q) E(P_y b) S=
  P_y b S,
  $$
  where the last step follows from the fact that
  $P_y \a(Q) E(P_y b) = P_y b$, which the reader may once again prove
by direct computation.
  It follows that
  $$
  \sum_{y\in\Gen} \Mat(x,y) S_yS_y^* bS =
  \sum_{y\in\Gen} \Mat(x,y) P_y bS =
  Q_x bS.
  $$
  We then conclude that the pair
  $$
  \(Q_x, \sum_{y\in\Gen} \Mat(x,y) S_yS_y^*\)
  $$
  is a redundancy and, observing that $\a(1)=1$ and hence the ideal
generated by $\R$ is the whole of $C(\X)$, we have that
  $$
  S_x^*S_x =
  Q_x \equiv
  \sum_{y\in\Gen} \Mat(x,y) S_yS_y^*
  $$
  modulo $\I$, as claimed.  It follows from the universal property of
$\OA$ \cite{\AnHR} that there exists a unique *-homomorphism
  $$
  \phi : \OA \to C(\X)\cp
  $$
  sending each canonical partial isometry generating $\OA$, which we
denote by $s_x$, to the class of $S_x$ in $C(\X)\cp$.

We now set out to define a map in the reverse direction.  In order to
do this recall from \cite{\QR} (see also \cite{\amena} and
\cite{\infinoa}) that $\OA$ is graded over the free group $\F =
\F(\Gen)$, and that the fiber over the identity group element may be
naturally identified with $C(\X)$.  Making this identification we have
that
  $s_x^*s_x = Q_x$ and
  $s_xs_x^* = P_x$.
  Defining
  $$
  s = \a(Q)^{-1/2}\sum_{x\in\Gen} s_x,
  $$
  we claim that $s$ satisfies \lcite{\DefineToeplitzExtension.i--ii}
with respect to $C(\X)$.  Speaking of
\lcite{\DefineToeplitzExtension.i} one may use the partial
crossed-product \cite{\infinoa} structure of $\OA$ to show that
  $
  s_x f = \a(f) s_x,
  $
  for all $x\in\Gen$ and $f\in C(\X)$ so that also
  $sf = \a(f) s$.  In order to verify
\lcite{\DefineToeplitzExtension.ii} observe that for $f\in C(\X)$ we
have
  $$
  s^* f s =
  \sum_{x,y\in\Gen} s_x^*\a(Q)\inv f s_y =
  \sum_{x\in\Gen} s_x^*\a(Q)\inv f s_x =
  \Tr(f),
  $$
  where the last step may also be proven using the partial
crossed-product structure of $\OA$.  By the universal property of
$\ATCP(C(\X),\a,\Tr)$ there exists a unique *-homomorphism
  $$
  \psi : \ATCP(C(\X),\a,\Tr) \to \OA
  $$
  which is the identity on $C(\X)$ and satisfies $\psi(S)=s$.

We claim that $\psi$ vanishes on the ideal $\I$ of
\lcite{\CrossedProductDef}.  In fact let $(a,k)$ be a redundancy and
notice that since $k\in \overline{A SS^*A}$ (from now on we will
denote $C(\X)$ simply by $A$), one has that $\psi(k)\in \overline{A
ss^*A}$.
  Observe however that
  $$
  Ass^*A \subseteq
  \sum_{x,y\in\Gen} As_xs_y^*A \subseteq
  \sum_{x,y\in\Gen} As_xs_y^*.
  $$
  The subspace $As_xs_y^*$ happens to be precisely the fiber over
$xy\inv$ when $x\neq y$, while $\bigoplus_{x\in\Gen} As_xs_x^* =
A$. It follows that $\sum_{x,y\in\Gen} As_xs_y^*$ is closed and hence
$\psi(k)$ must have the form
  $$
  \psi(k) = \sum_{x,y\in\Gen} k_{xy}s_x s_y^*,
  $$
  where the coefficients $k_{xy}$ are in $A$.  By hypothesis we have
that $abS = kbS$ and hence
  $$
  abs = \psi(k)bs,
  $$
  for all $b\in A$, which translates into
  $$
  a b \a(Q)^{-1/2}\sum_{x\in\Gen} s_x =
  \sum_{x,y,z\in\Gen} k_{xy}s_x s_y^*
  b
  \a(Q)^{-1/2}s_z.
  $$
  Substituting $b$ for $b \a(Q)^{-1/2}$,
  observing that $s_y^* b s_z = 0$ when $y\neq z$, and
  projecting on the fiber over $x$, we have that
  $$
  a b s_x =
  \sum_{y\in\Gen} k_{xy}s_x s_y^* b s_y,
  \eqno{(\dagger)}
  $$
  for all $x\in \Gen$.  Taking $b=P_y$, with $y\in\Gen$, we get
  $$
  a P_y s_x =
  k_{xy}s_x s_y^* P_y s_y =
  k_{xy}s_x Q_y.
  $$
  When $x\neq y$ the left hand side vanishes and hence
  $$
  k_{xy}s_x s_y^*=
  k_{xy}s_x Q_ys_y^*=
  0.
  $$
  It follows that
  $$
  \psi(k) =
  \sum_{x\in\Gen} k_{xx}s_x s_x^* =
  \sum_{x\in\Gen} k_{xx}P_x,
  $$
  and also that \lcite{$\dagger$} reduces to
  $$
  a b s_x =
  k_{xx}s_x s_x^* b s_x.
  $$
  Multiplying this on the right by $s_x^*$ leads to
  $$
  a bP_x =
  k_{xx}s_x s_x^* b s_x s_x^* =
  k_{xx} b P_x.
  $$
  With $b=1$ we get
  $$
  a =
  \sum_{x\in\Gen} aP_x =
  \sum_{x\in\Gen} k_{xx}P_x =
  \psi(k).
  $$
  Therefore $\psi(a-k)=0$ and hence $\psi$ vanishes on $\I$ as claimed
and hence factors through the quotient yielding a *-homomorphism
  $$
  \widetilde\psi : C(\X)\cp \to \OA
  $$
  which may now be proven to be the inverse of $\phi$.  This concludes
the proof.
  \proofend

It should be remarked that in the last part of the proof above we have
\stress{proved} that $\psi(a-k)=0$ for \stress{all} redundancies.
This is perhaps an indication that one should in fact mod out all
redundancies in Definition \lcite{\CrossedProductDef} when $\a(1)=1$.

\endgroup

\end
\references

\bibitem{\Baladi}
  {V. Baladi}
  {Positive transfer operators and decay of correlations}
  {Advanced Series in Nonlinear Dynamics vol.~16, World Scientific,
2000}

\bibitem{\Blackadar}
  {B. Blackadar}
  {Shape theory for $C^*$-algebras}
  {\sl Math. Scand. \bf 56 \rm (1985), 249--275}

\bibitem{\Cuntz}
  {J. Cuntz}
  {Simple $C^*$-algebras generated by isometries}
  {\sl Comm. Math. Phys. \bf 57 \rm (1977), 173--185}

\bibitem{\CuntzTwo}
  {J. Cuntz}
  {The internal structure of simple C*-algebras}
  {\sl Operator algebras and applications, Proc. Symp. Pure Math. \bf
38 \rm (1982), 85-115}

\bibitem{\CK}
  {J. Cuntz and W. Krieger}
  {A Class of C*-algebras and Topological Markov Chains}
  {\sl Inventiones Math. \bf 56 \rm (1980), 251--268}

\bibitem{\Deaconu}
  {V. Deaconu}
  {Groupoids associated with endomorphisms}
  {\sl Trans. Amer. Math. Soc. \bf 347 \rm (1995), 1779--1786}

\bibitem{\DeaconuMuhly}
  {V. Deaconu and P. Muhly}
  {C*-algebras associated with Branched Coverings}
  {preprint, University of Nevada. To appear in {\sl
Proc. Amer. Math. Soc}}

\bibitem{\DR}
  {S. Doplicher and J. E. Roberts}
  {Endomorphisms of C*-algebras, cross products and duality for
compact groups}
  {\sl Ann. Math. \bf 130 \rm (1989), 75--119}

\bibitem{\amena}
  {R. Exel}
  {Amenability for {F}ell Bundles}
  {\sl J. reine angew. Math. \bf 492 \rm (1997),
41--73. [funct-an/9604009]}

\bibitem{\infinoa}
  {R. Exel and M. Laca}
  {{C}untz--{K}rieger Algebras for Infinite Matrices}
  {\sl J. reine angew. Math. \bf 512 \rm (1999),
119--172. [funct-an/9712008]}

\bibitem{\AnHR}
  {A. an Huef and I. Raeburn}
  {The ideal structure of Cuntz-Krieger algebras}
  {\sl Ergodic Theory Dyn. Syst. \bf 17 \rm (1997), 611--624}

\bibitem{\Murphy}
  {G. J. Murphy}
  {Crossed products of C*-algebras by endomorphisms}
  {\sl Integral Equations Oper. Theory \bf 24 \rm (1996), 298--319}

\bibitem{\Paschke}
  {W. L. Paschke}
  {The Crossed Product of a $C^*$-algebra by an Endomorphism}
  {\sl Proc. Amer. Math. Soc. \bf 80 \rm (1980), 113--118}

\bibitem{\Pim}
  {M. V. Pimsner}
  {A class of C*-algebras generalizing both Cuntz-Krieger algebras and
crossed products by $\bf Z$}
  {\sl Fields Inst. Commun. \bf 12 \rm (1997), 189--212}

\bibitem{\PV}
  {M. Pimsner and D. Voiculescu}
  {Exact sequences for $K$-groups and Ext-groups of certain
cross-products $C^*$-algebras}
  {\sl J. Oper. Theory \bf 4 \rm (1980), 93--118}

\bibitem{\QR}
  {J. Quigg and I. Raeburn}
  {Characterisations of crossed products by partial actions}
  {\sl J. Oper. Theory \bf 37 \rm (1997), 311--340}

\bibitem{\Renault}
  {J. Renault}
  {A groupoid approach to $C^*$-algebras}
  {Lecture Notes in Mathematics vol.~793, Springer, 1980}

\bibitem{\Ruelle}
  {D. Ruelle}
  {The thermodynamic formalism for expanding maps}
  {\sl Commun. Math. Phys. \bf 125 \rm (1989), 239--262}

\bibitem{\Stacey}
  {P. J. Stacey}
  {Crossed products of C*-algebras by *-endomorphisms}
  {\sl J. Aust. Math. Soc., Ser. A \bf 54 \rm (1993), 204--212}
